\documentclass[11pt]{article}

\usepackage{amsmath,amssymb,bbm,pslatex,xcolor,graphics}

\usepackage[dvips]{graphicx}
\usepackage{times,rotate}

          \newtheorem{theorem}{Theorem}[section]
          
          \newtheorem{lemma}[theorem]{Lemma}
          
          \newtheorem{definition}{Definition}[section]
          \newtheorem{corollary}[theorem]{Corollary}

\catcode`\<=\active \def<{
\fontencoding{T1}\selectfont\symbol{60}\fontencoding{\encodingdefault}}
\catcode`\>=\active \def>{
\fontencoding{T1}\selectfont\symbol{62}\fontencoding{\encodingdefault}}
\catcode`\|=\active \def|{
\fontencoding{T1}\selectfont\symbol{124}\fontencoding{\encodingdefault}}

\newcommand{\tmop}[1]{\ensuremath{\operatorname{#1}}}

\def\ed{\end{document}}

\def\up{\underline{P}}

\def\cal{\mathcal}
\def\widecheck{\check}

\def\qed{\mbox{}\hfill Q.E.D.}
\def\game{\Game }
\def\Df{\mathop{Def}}

\begin{document}

\catcode`\<=\active \def<{
\fontencoding{T1}\selectfont\symbol{60}\fontencoding{\encodingdefault}}
\catcode`\>=\active \def>{
\fontencoding{T1}\selectfont\symbol{62}\fontencoding{\encodingdefault}}
\catcode`\|=\active \def|{
\fontencoding{T1}\selectfont\symbol{124}\fontencoding{\encodingdefault}}
\newcommand{\tmrsub}[1]{\ensuremath{_{\textrm{#1}}}}
\newcommand{\tmstrong}[1]{\textbf{#1}}


\newcommand{\hi}{version Sep 23 15}
\newcommand{\mb}{\mbox{ }}
\newcommand{\abo}[2]{\ensuremath{#1 <_{\mathcal{O}}   #2}}
\newcommand{\calo}{ {\mathcal{O}}}
\newcommand{\calm}{\ensuremath{\mathcal{M}} }
\newcommand{\cali}{\ensuremath{\mathcal{I}} }
\newcommand{\calj}{{\mathcal{J}}}
\newcommand{\calk}{ {\mathcal{K}}}
\newcommand{\calc}{C}
\newcommand{\gott}{\ensuremath{\mathfrak{T}}}
\newcommand{\GRP}{\ensuremath{(\ensuremath{\operatorname{GRP}}_{0} )}}
\newcommand{\LST}{\ensuremath{\call_{\dot{\in}}}}
\newcommand{\ent}{{{\em Entscheidungsproblem\/}}}
\newcommand{\bara}{\bar{{\alpha}}}
\newcommand{\bp}{\bar{{\pi}}}
\newcommand{\ts}{\ensuremath{T^{\ast}}}
\newcommand{\tsptwo}{\left.{\emp}\right) \
\ensuremath{_{\textrm{p\tmrsub{2}}}}}
\newcommand{\tptwo}{\ensuremath{T^{\ast} ( \varnothing )_{p_{2}}}}
\newcommand{\pe}{\ensuremath{P^{\ensuremath{\operatorname{\mathsf{eJ}}}}_{e}}}
\newcommand{\tmt}[1]{\ensuremath{\mathbbm{T}^{#1}}}
\newcommand{\p}[1]{\ensuremath{P^{\ensuremath{\operatorname{\mathsf{eJ}}}}_{#1}}}
\newcommand{\zx}[1]{\ensuremath{^{#1} \zeta}}
\newcommand{\sx}[1]{\ensuremath{^{#1} \Sigma}}
\newcommand{\tx}[1]{\ensuremath{^{#1} T}}
\newcommand{\cst}{\ensuremath{\macro{C^{\ast}}}}
\newcommand{\tmcr}{\ensuremath{\curvearrowright}}
\newcommand{\cl}{\ensuremath{\curvearrowleft}}
\newcommand{\ou}{{\"o}}
\newcommand{\au}{{\"a}}
\newcommand{\uu}{{\"u}}
\newcommand{\go}{G{\"o}del}
\newcommand{\phiav}[1]{\ensuremath{\Phi^{\alpha}_{#1}}}
\newcommand{\phiaf}{\ensuremath{\Phi_{\ensuremath{\operatorname{fp}}}^{\alpha}}}
\newcommand{\phiabp}[2]{\ensuremath{\Phi_{\beta +1}^{(  \alpha )}}}
\newcommand{\phiab}[2]{\ensuremath{\Phi_{#2}^{#1}}}
\newcommand{\phivf}[1]{\ensuremath{\Phi_{\ensuremath{\operatorname{fv}}}^{#1}}}
\newcommand{\rar}{{\rightarrow}}
\newcommand{\midd}{ {\mid} }
\newcommand{\con}{{\smallfrown}}
\newcommand{\lt}{{\leq}\ensuremath{_{\textrm{T}}}}
\newcommand{\qu}{Q}
\newcommand{\pb}{\ensuremath{\Psi_{\beta}}}
\newcommand{\oy}{<\ensuremath{_{\textrm{y\tmrsub{0}}}}}
\newcommand{\oyr}{\ensuremath{>_{y_{0}}}}
\newcommand{\oye}{{\leq}\ensuremath{_{\textrm{y\tmrsub{0}}}}}
\newcommand{\ptwo}{{{\em II\/}}}
\newcommand{\pone}{{{\em I\/}}}
\newcommand{\uhr}{{\upharpoonright}}
\newcommand{\2}{{{\em II\/}}}
\newcommand{\1}{{{\em I\/}}}
\newcommand{\wff}{\ensuremath{\ensuremath{\operatorname{wff}}(y_{0} )}}
\newcommand{\nod}{{\noindent}}
\newcommand{\diam}{\ensuremath{}{\diamond}}
\newcommand{\lalp}{ {\pa{L\ensuremath{_{\textrm{{\alpha}}}},{\in}}}}
\newcommand{\cit}{{\macro{}{{\cite{\arg{k}}}}}}
\newcommand{\ul}{\ensuremath{\ulcorner}}
\newcommand{\ur}{\ensuremath{\urcorner}}
\newcommand{\cs}{\ensuremath{\ul \sigma \ur}}
\newcommand{\ct}{\ensuremath{\ul \tau \ur}}
\newcommand{\cso}{\ensuremath{\ul \sigma_{0} \ur}}
\newcommand{\co}[1]{\ensuremath{\ulcorner #1 \urcorner}}
\newcommand{\gp}[1]{\ensuremath{\prec #1 \succ}}
\newcommand{\tmfi}{\ensuremath{\ensuremath{\operatorname{Field}}( \prec )}}
\newcommand{\al}{\ensuremath{\aleph}}
\newcommand{\alo}{\ensuremath{\aleph_{\omega}}}
\newcommand{\cub}{c.u.b.}
\newcommand{\dff}{\ensuremath{\ensuremath{\operatorname{Def}}}}
\newcommand{\bu}{\ensuremath{\bullet}}
\newcommand{\nat}{\ensuremath{\mathbbm{N}}}
\newcommand{\re}{\ensuremath{\mathbbm{R}}}
\newcommand{\baire}{\ensuremath{\mathbbm{N}}\textsuperscript{\ensuremath{\mathbbm{N}}}}
\newcommand{\bai}{\textsuperscript{\ensuremath{\omega}}\ensuremath{\omega}}
\newcommand{\linf}{{\leq}\ensuremath{_{\textrm{{\infty}}}}}
\newcommand{\linfs}{<\ensuremath{_{\textrm{\text{{\tmstrong{\ensuremath{\infty}}}}}}}}
\newcommand{\leinf}{{\leq}\ensuremath{_{\textrm{\text{{\tmstrong{\ensuremath{\infty}}}}}}}}
\newcommand{\cant}{2\textsuperscript{\ensuremath{\mathbbm{N}}}}
\newcommand{\can}{\textsuperscript{\ensuremath{\omega}}2}
\newcommand{\rest}{{\upharpoonright}}
\newcommand{\emp}{{\varnothing}}
\newcommand{\pa}[1]{\ensuremath{\langle #1 \rangle}}
\newcommand{\wst}{\ensuremath{S}}
\newcommand{\ws}[2]{\ensuremath{S}\textsuperscript{#1}\ensuremath{_{\textrm{#2}}}}
\newcommand{\wsone}[1]{\ensuremath{S}\textsuperscript{1}\ensuremath{_{\textrm{#1}}}}
\newcommand{\wstwo}[1]{\ensuremath{S}\textsuperscript{2}\ensuremath{_{\textrm{#1}}}}
\newcommand{\wthree}[1]{\ensuremath{S}\textsuperscript{3}\ensuremath{_{\textrm{#1}}}}
\newcommand{\es}[2]{\ensuremath{E}\textsuperscript{#1}\ensuremath{_{\textrm{#2}}}}
\newcommand{\esone}[1]{\ensuremath{E}\textsuperscript{1}\ensuremath{_{\textrm{#1}}}}
\newcommand{\estwo}[1]{\ensuremath{E}\textsuperscript{2}\ensuremath{_{\textrm{#1}}}}
\newcommand{\wet}{\ensuremath{E}}
\newcommand{\games}{\ensuremath{\Game \Sigma^{0}_{3}}}
\newcommand{\szero}{\ensuremath{\Sigma^{0}_{3}}}
\newcommand{\hright}{{\hookrightarrow}}
\newcommand{\call}{L}
\newcommand{\callt}{L\ensuremath{_{\textrm{T}}}}
\newcommand{\calt}{T}
\newcommand{\cald}{D}
\newcommand{\calu}{\ensuremath{\mathcal{U}}}
\newcommand{\lst}{\ensuremath{\mathcal{L}_{\dot{\in}}}}
\newcommand{\calw}{W}
\newcommand{\calv}{V}
\newcommand{\vp}{\ensuremath{\varphi}}
\newcommand{\om}[1]{\ensuremath{\omega_{1} ( #1 ) }}
\newcommand{\da}{{\downarrow}}
\newcommand{\ua}{{\uparrow}}
\newcommand{\la}{{\langle}}
\newcommand{\ra}{{\rangle}}
\newcommand{\ran}{\text{ran}}
\newcommand{\dom}{\text{dom}}
\newcommand{\back}{{\backslash}}
\newcommand{\power}{P}
\newcommand{\mar}{{\marginpar{\ensuremath{\rightarrow}}}}
\newcommand{\bm}{{\boldmath}}
\newcommand{\opistol}{\text{\ensuremath{O^{\pistol}}}}
\newcommand{\ie}{{\itshape{i.e.}}{\hspace{0.25em}}}
\newcommand{\eg}{{\itshape{e.g.}}{\hspace{0.25em}}}
\newcommand{\etc}{{\itshape{etc.}}{\hspace{0.25em}}}
\newcommand{\via}{{\itshape{via }}{\hspace{0.25em}}}
\newcommand{\cf}{{\itshape{cf. }}{\hspace{0.25em}}}
\newcommand{\mm}{{\itshape{m.m.}}{\hspace{0.25em}}}
\newcommand{\an}{{\wedge}}
\newcommand{\bigconj}{\ensuremath{\bigwedge \hspace{-1.0em} \bigwedge}}
\newcommand{\tmk}{\text{``}}
\newcommand{\proves}{{\vdash}}
\newcommand{\all}{{\forall}}
\newcommand{\Equi}{{\Longleftrightarrow}}
\newcommand{\equi}{{\longleftrightarrow}}
\newcommand{\ex}{{\exists}}
\newcommand{\Imp}{{\Rightarrow}}
\newcommand{\sset}{\hspace{.11em}\subseteq \hspace{.11em}}
\newcommand{\lr}{{\leftrightarrow}}
\newcommand{\imp}{{\longrightarrow}}
\newcommand{\ri}{{\rightarrow}}
\newcommand{\nequiv}{{\not\equiv}}
\newcommand{\wt}{\tilde{\}}}
\newcommand{\finold}{\ensuremath{[ \alpha ]^{< \omega}}}
\newcommand{\rem}{{\noindent}{\bfseries{Remark: }}}
\newcommand{\fin}[1]{\ensuremath{[  ]^{< \omega}}}
\newcommand{\dlimtwo}{\ensuremath{\raisebox{-1.25\ensuremath{\operatorname{ex}}}{\overset{Lim}{\longrightarrow}
\hspace{0.25em}}}}
\newcommand{\dlim}{\text{{\raisebox{-1.25ex}{\ensuremath{\overset{Lim}{\longrightarrow}
\hspace{0.25em}}}}}}
\newcommand{\qedtwo}[1]{QED(#1)}
\newcommand{\pf}{{\noindent}{\textbf{Proof: }}}
\newcommand{\dfs}{\ensuremath{=_{\ensuremath{\operatorname{df}}}}}
\newcommand{\edfs}{\ensuremath{  \Equi_{\ensuremath{\operatorname{df}}}  }}
\newcommand{\pr}{\ensuremath{\prec}}
\newcommand{\pre}{\ensuremath{\preceq}}
\newcommand{\js}{\ensuremath{J_{s}}}
\newcommand{\jsb}{\ensuremath{J_{\bar{s}}}}
\newcommand{\jst}[3]{ \ensuremath{\langle J^{#1}_{#2} , #3 \rangle}}
\newcommand{\jnb}{J\ensuremath{_{\textrm{\bar{{\nu}}}}}}
\newcommand{\jn}{J\ensuremath{_{\textrm{{\nu}}}}}
\newcommand{\lam}[2]{\ensuremath{\Lambda ( #1 , #2 )}}
\newcommand{\nb}{\ensuremath{\bar{\nu}}}
\newcommand{\nbn}{\ensuremath{\bar{\nu} \prec \nu}}
\newcommand{\ns}{{\nu}\ensuremath{_{\textrm{s}}}}
\newcommand{\tmsb}{\ensuremath{\bar{s}}}
\newcommand{\nsb}{\ensuremath{\nu_{\bar{s}}}}
\newcommand{\maps}{\ensuremath{f: \bar{s} \Longrightarrow s}}
\newcommand{\mapnr}{\ensuremath{f \Longrightarrow \nu}}
\newcommand{\mapn}{\ensuremath{f: \bar{\nu} \Longrightarrow \nu}}
\newcommand{\mapg}[3]{\ensuremath{#1 : \overline{#2} \Longrightarrow #3}}
\newcommand{\onektuple}[1]{\ensuremath{#1_{1} , \ldots ,  #1_{k}}}
\newcommand{\zeroktuple}[1]{\ensuremath{#1_{0} , \ldots  , #1_{k}}}
\newcommand{\onentuple}[1]{\ensuremath{#1_{1} , \ldots  , #1_{n}}}
\newcommand{\zerontuple}[1]{\ensuremath{#1_{0} , \ldots  , #1_{n}}}
\newcommand{\ntuple}[2]{\ensuremath{#1_{1} , \ldots  , #1_{#2}}}
\newcommand{\otuple}[2]{\ensuremath{#1_{0} , \ldots  , #1_{#2}}}
\newcommand{\signk}{\ensuremath{\Sigma_{k}^{(n)}}}
\newcommand{\signo}{\ensuremath{\Sigma_{1}^{(n)}}}
\newcommand{\bsigno}{\ensuremath{\ensuremath{\boldsymbol{\Sigma}}_{1}^{(n)}}}
\newcommand{\bsigmo}{\ensuremath{\ensuremath{\boldsymbol{\Sigma}}_{1}^{(m)}}}
\newcommand{\sigmk}[1]{\ensuremath{\Sigma_{#1}^{(n)}}}
\newcommand{\sigmnk}[2]{\ensuremath{\Sigma_{#2}^{( #1 )}}}
\newcommand{\sigms}{ \ensuremath{\Sigma^{\ast}}}
\newcommand{\bsigms}{{\textbf{ \ensuremath{\Sigma^{\ast}}}}}
\newcommand{\lnn}{ \ensuremath{\Lambda (q, \nu )}}
\newcommand{\lns}{ \ensuremath{\Lambda (q,s)}}
\newcommand{\lnsp}{ \ensuremath{\Lambda^{+} (q,s)}}
\newcommand{\lnsl}{ \ensuremath{\Lambda (q,s| \lambda )}}
\newcommand{\fb}[1]{\ensuremath{\lambda (f_{( #1 ,q, \nu )} )}}
\newcommand{\fbb}[2]{\ensuremath{\lambda (f_{( #1 ,q,  #2 )} )}}
\newcommand{\fbs}[1]{\ensuremath{\lambda (f_{( #1 ,q,s)} )}}
\newcommand{\lis}{\ensuremath{l^{i}_{s}}}
\newcommand{\lhs}{\ensuremath{l^{i}_{\eta s}}}
\newcommand{\lls}{\ensuremath{l^{i}_{\lambda s}}}
\newcommand{\llf}{\ensuremath{\lambda = \lambda (f)}}
\newcommand{\csp}{C\ensuremath{_{\textrm{s}}}\textsuperscript{+}}
\newcommand{\csl}{C\ensuremath{_{\textrm{s|{\lambda}}}}}
\newcommand{\rn}{ \ensuremath{\rho^{n+1}_{M} > \kappa}}
\newcommand{\rhn}[1]{{\rho}\textsuperscript{#1}\ensuremath{_{\textrm{M}}}}
\newcommand{\rhnm}[2]{{\rho}\textsuperscript{#1}\ensuremath{_{\textrm{#2}}}}
\newcommand{\eqn}{{\omega}{\rho}\textsuperscript{n+1}\ensuremath{_{\textrm{M}}}{\leq}{\kappa}<{\omega}{\rho}\textsuperscript{n}\ensuremath{_{\textrm{M}}}}
\newcommand{\hnm}{H\textsuperscript{n}\ensuremath{_{\textrm{M}}}}
\newcommand{\oddpagetext}[1]{\newcommand{\pageoddheader}{{\small }}}
\newcommand{\evenpagetext}[1]{\newcommand{\pageevenheader}{{\small }}}


\title{
  {The Ramified Analytical Hierarchy using Extended Logics}
}

\author{  P.D. Welch\\  \small {School of Mathematics,} \\ \small University of Bristol,  \\ \small Bristol, BS8 1TW, England}

\maketitle 

\begin{abstract}  
  
The use of Extended Logics to replace ordinary second order definability in Kleene's {\em Ramified Analytical Hierarchy} is investigated.  This mirrors a similar investigation of Kennedy, Magidor and V\"a\"an\"anen \cite{KeMaVa2016} where G\"odel's universe $L$ of constructible sets is subjected to similar variance. Enhancing second order definability allows models to be defined which may or may not coincide with the original Kleene hierarchy in domain. Extending the logic with game quantifiers, and assuming strong axioms of infinity, we obtain {\em minimal correct} models of analysis.   A wide spectrum of models can be so generated from abstract definability notions: one may  take an abstract   Spector Class and  extract an extended logic for it.
The resultant structure is then a minimal model of the given kind of definability.

\end{abstract}

\maketitle 

\section{Introduction}

This paper arises out of questions of Kennedy {\cite{Ke2013}} and the authors of {\cite{KeMaVa2016}}.\footnote{Keywords: definability, analytical hierarchy, determinacy; Classifications: 03E60, 03E15, 03E47.}  From the first paper cited:

\begin{quote}
We read G\"odel's 1946 lecture as an important but perhaps overlooked step in this line of thought [{\em concerning formalism freeness}], not with respect to language necessarily  \mbox{[ . . .  ]} but with respect to formalization altogether; in particular we will interpret G\"odel there as making the suggestion, albeit in a preliminary form, that tests of robustness analogous to that which is implicit in the Church-Turing Thesis be developed, not for the notion of computable function but for the concept of definability - witnessing its formalism independence, as it were.
\end{quote}

The second paper cited looks at building inner model hierarchies as G\"odel did, but instead using definability of the models using languages with extended quantifiers, rather than just first order logic, thereby testing the extent to which $L$ was indeed ``formalism-free'' or independent of the logic used. Of course in some cases the altered formalism did indeed return $L$, however in  others this was not the case, with new and interesting inner models arising.  We seek to follow up the question raised by the following suggestion (\cite{Ke2014}, also raised in the previously cited paper).

\begin{quote}The method can be implemented not just for definability in the sense of $L$  or $HOD$ as was done in \cite{KeMaVa2016}, 
but in other settings as well. Thinking beyond definability toward other canonical concepts, one might also consider this varying the underlying logic also in other contexts. In fact any logical hierarchy, e.g., Kleene's ramified hierarchy of reals is amenable to this treatment, conceivably. Suitable notions of confluence and grounding must be formulated on a case by case basis.\end{quote}

We shall attempt then to give one answer to a question of whether something similar might be done for the {\em Ramified Analytical Hierarchy} of Kleene. We shall first illustrate the methodology for inner models from \cite{KeMaVa2016}, immediately followed by an example: the ``{\em cof-$\omega$}'' quantifier, $\mathsf{Q}^{cf}_{\omega}$.  We then give the definition of the Ramified Analytical Hierarchy (and note why the  $\mathsf{Q}^{cf}_{\omega}$ is inappropriate in this setting).  We proceed then to talk about generalized quantifiers that are appropriate. This theory goes back to Aczel, \cite{Ac75}, and to Moschovakis (see for example \cite{M}), and we focus on examples given by {\em game quantifiers}. As we shall see these provide a rich source of differing hierarchies.

  In Section 2 we define the {\em minimum correct model of analysis}. This uses the game quantifier and a background theory of Projective Determinacy. It is quite possible to define a correct model using a different background theory such as is done in \cite{Sh72} assuming `$V=L$'.  It is a feature of  \cite{KeMaVa2016} that varying the background theories, such as taking $V=L$ or $L[\mu]$, or forcing extensions, or, ... gives rise to different versions of a model defined by the same logic. It may be argued that only when we have sufficiently strong axioms can we obtain the definitive model.  The same is true for correct models of ramified analysis. 

In Section 4 we parallel the result of Gandy and Putnam (\cite{BHP66}) that found the least level of an inner model (namely G\"odel's $L$) whose reals corresponded to the least $\beta$-model of analysis. We identify the level, $Q$, of the least iterable inner model with infinitely many Woodin cardinals whose reals correspond to those of $\up^{Proj}$, the minimal fully correct model of analysis (`fully correct' in the sense that each $\Pi^{1}_{n}$ formula is absolute between the model $\up^{Proj}$ and standard model of analysis).

Thus:

\setcounter{section}{4}\setcounter{theorem}{7}
\begin{theorem}   $P^{Proj}=\re\cap Q$.
\end{theorem}  
 In the final section we argue that, using a representation theorem of Harrington, for almost any notion of `definability' for sets of integers in some general sense, we can find an extended logic for it - this is the contents of the following theorem: 
\setcounter{section}{6}\setcounter{theorem}{1}
\begin{theorem}\label{6.2} Let $\Gamma\sset \power(\nat)$ be a Spector class, with corresponding quantifier $\mathsf{Q}= \mathsf{Q}_{\Gamma}$ from  Harrington's theorem. Then there is a minimum model of analysis $\up^{\Gamma}$ which is closed under positive inductions in $\call^{\mathsf{Q}}$, and so that for any $X\in P^{\Gamma}$ we have that $\Gamma(X)$ (the Spector class relativised to $X$) is contained in $\up^{\Gamma}$.
\end{theorem}
\setcounter{section}{1}\setcounter{theorem}{1}
 This might be considered a possible maximal answer to Kennedy's question. We close with some further open questions. 

\section{Extended logics for set theory and for ramified analysis}

We first give an example from inner model theory.  This is done by varying the G\"odelian $\mathop{Df}$ operator in the recursive definition of the constructible universe. By `a logic' is meant a set of sentences $S$ in a {\em language }  and a {\em truth}  (and so a satisfaction) predicate $T$ for them.  It is intended that $T$ give full information for any structure (of the appropriate signature) $\calm$ as to whether $\varphi$ holds in $\calm$ or not, thus we may take $T$ as a function $T(\calm,\varphi)$ with values in $2$. The {\em logic} $\call^{\ast}$ is then the pair $(S,T)$. (We may slip up  in the sequel and speak just of  $\call^{\ast}$ as a language, but the reader will know what is meant.) For ordinary first order logic the following is just G\"odel's definition of his $L$.

\begin{definition} If $M$ is a set, let $\Df _{\call^{\ast}} (M)$ 
denote the set of all sets of the form $$X = \{a \in M \mid (M,\in) \models \varphi(a,b)\}$$ where $\varphi(x,y)$ is an arbitrary formula of the logic $\call^{\ast}$ and $b \in M$. We define a hierarchy 
$(L'_{\alpha})$ of sets constructible using $\call^{\ast}$ as follows:
\begin{center}
$L'_{0}=\emp$\,;\\
$L'_{\nu} = \bigcup_{\alpha < \nu} L'_{\alpha}$ for $Lim(\nu)$;\\
$L'_{\alpha+1}= \Df_{\call^{\ast}}(L'_{\alpha}) $ \,and $L' = \bigcup_{\alpha \in On} L'_{\alpha}.$\\
\end{center}
\end{definition}

 One example: $\call^{\ast}$ is the language  $L_{\dot \in,\dot =}$ of set theory together with a {\em cofinality $\omega$ quantifier}.\\

$ (M,\in)\models \mathsf{Q}^{cf}_{\omega}xy\varphi(x,y,a) \Leftrightarrow $\\
\mb \hfill$ \{(c,d) \mid (M ,\in)\models \varphi (c,d,a)\}
\mbox{ is a linear order of cofinality } \omega.$\\

 The model built here, $L'$, can be shown to be precisely $L[C^{\omega}]$ where $C^{\omega}$ is the class of ordinals of cofinality $\omega$.  One should note that  then $L'$ `knows' that certain ordinals have cofinality $\omega$ in $V$, but they need not have countable cofinality in $L'$: the construction does not provide a cofinal $\omega$ sequence for each $\alpha\in C^{\omega}$. 

 It is thus important to remark that this is not an absolute notion: the right hand side is evaluated in $V$, not in any sense in the final, as yet to be built, model $L'$.  New information is thus imported from $V$ into the construction of $L'$.  
 
 To state the obvious: this importation of new information is an essential feature of each of the logics/quantifiers of {\cite{KeMaVa2016}} that build something different from $L$.
 
We take up the question of the second quotation above concerning analogous constructions   for the Ramified Analytical Hierarchy.
This is expressed in a suitable language $\call^{2}$ for analysis. (By `analysis' we mean an axiomatisation of second order number theory, such as $\mathsf{Z_{2}}$, see for example \cite{Si99}.) Such a language is appropriate for structures of the form
$$\calm = (M, \nat, +, \times, 0,', \ldots)$$
where $M\sset \power({\nat})$.  Thus $\call^{2}$ contains {\em number variables} $x,y,z,\ \ldots$ and second order {\em set variables} $X,Y,Z, \ldots$, for sets of numbers, and quantifiers of both kinds: $\ex
x, \all y, \ldots \ex X, \all Y \ldots$, as well as function symbols for $0, +, \times, \ldots $, \etc We let  be the full model of analysis be $\cal{Z}= (\power(\nat), \nat, +, \times, 0,', \ldots)$. We also identify $\re$ with $\power(\nat)$ and use both interchangeably.

\begin{definition}[Kleene]{\em(\cite{Kl59})}  Define by recursion on $\alpha$, $P_{\alpha}\sset \power(\nat)$: 
\begin{center}
$P_{0}=\emp$;\\
$P_{\lambda} = \bigcup_{\alpha <\lambda}P_{\alpha} $ for $Lim(\lambda)$;\\

$P_{\alpha +1}= \{Y\sset \nat \mid Y \mbox{ is definable in } \call^{2} \mbox{ over }   \underline{P_{\alpha}} = (P_{\alpha}, \nat, +,\times , 0, ', \ldots )\}$
\end{center}            

\end{definition}
On cardinality grounds there will be a {\em fixed point} ${P} =_{df} {P}_{\beta_{0}} = {P}_{\beta_{0}+1}$ \\

It was Cohen in \cite{Coh63} who first showed that this ordinal $\beta_{0}$ was countable. (How could this have been ever in question given the L\"owenheim-Skolem theorem? Presumably in this relatively early period people were simply unused to working with, and deploying arguments about, the absoluteness of various constructions to the constructible hierarchy. Using this absoluteness, taking a countable substructure of a part of the universe containing the hierarchy up to the fixed point would have revealed the latter's countability.)
Cohen also conjectured that $\underline{P}_{\beta_{0}}$ formed the minimal $\beta$-model of analysis.  (A model of a fragment of second order theory is a $\beta${\em -model} if it is ``correct'' or absolute for $\Pi^{1}_{1}$ expressions. For a discussion of  these notions  \cf\cite{Si99}.) This was subsequently proven independently by Gandy and by Putnam (\cite{BHP66}).

 Can we introduce non-standard quantifiers here and see what models, now not of {\em set theory}, but of {\em analysis} we can build?
 One problem, or difference, is that compared to the universe of sets, the Gandy-Putnam result shows that model $\underline{P}$ is tiny:
$$P_{\beta_{0}} = \power(\nat)\cap L_{\beta_{0}}$$ and $ L_{\beta_0}$ is the least level of the G\"odel hierarchy which is a $ZF^{-}$ model. 
(Later work by Boolos and Putnam \cite{BoPu68} gave a level by level analysis of this hierarchy and showed that $P_{\alpha} = \power(\nat)\cap L_{\alpha}$ for $\alpha\leq \beta_{0}$. From today's perspective these results seem again entirely straightforward, but we are standing on the shoulders of Jensen's magisterial fine structural analysis of the constructible universe.)

Hence $L_{\beta_{0}}\models
$``$V=HC$'' so there the $\mathsf{Q}^{cf}_{\omega}$ and other cardinality quantifiers are not going to get any traction in this region of analysis as everything is countable and provably so. Indeed for every level $\alpha <\beta_{0} ,\,L_{\alpha+1}\models $``$card(L_{\alpha})=\omega$''. So every ordinal is immediately collapsed and made cofinal with $\omega$.  It thus seems that set-theoretic based quantifiers are not going to be appropriate.
However, we may find in the papers of Aczel relevant analytically (so to speak) generalised quantifiers suitable for second order number theory. These are indeed entirely general, see \cite{Ac75}. Such a {\em quantifier} on $\power(\nat)$ is a set $\mathsf{Q}$ with $\emp\subsetneq \mathsf{Q} \subsetneq \power(\nat)$ which is {\em monotone}, that is $X\in \mathsf{Q} \wedge X\sset Y \imp Y\in \mathsf{Q}$. It is usual to write interchangeably
$$Y\in \mathsf{Q} \,\equi\, \mathsf{Q}(Y)\,\equi\, \mathsf{Q}x(Y(x)).$$
The {\em dual} of $\mathsf{Q}$ is the quantifier $\widecheck {\mathsf{Q}}$ given by:
$$\widecheck {\mathsf{Q}}xY(x)\, \equi \, \neg \mathsf{Q}x\neg Y(x).$$ 
Example (i) The Urexample of course is $\exists$ with dual $\forall$. \\
Example (ii) The {\em Suslin quantifier}. We assume $\la \,\ra:\bigcup_{n}\nat^{n}\equi\nat$ is a recursive bijection. Then this is defined by:
$${\cal{S}}uP(u) \Equi \all x_{0}\all x_{1}\all x_{2}\, \cdots \, \ex n P(\pa{x_{0},\ldots\, , x_{n}}).$$
\nod It's dual, $\widecheck {\cal{S}}$, is usually written $ {\cal{A}}$.  As can be seen, for an arithmetic $P$ this is equivalent to a $\Pi^{1}_{1}$ expression, and for $P$ $\Sigma^{1}_{n}$ it remains $\Sigma^{1}_{n}$ but only for $n\geq 2$.

A common quantifier is the {\em game quantifier}. An {\em infinite two person perfect information game} $G_{A}$ with $A\subseteq  \nat^{\nat}$ and players $I$, $I\! I$ is set up as usual:\\

$I \quad n_{0}\quad n_{1} \quad \ldots $

$I\! I \quad\quad m_{0}\quad m_{1}\quad \ldots\quad\quad \quad z=(n_{0},m_{0},n_{1},m_{1}, \ldots )$\\

$I$ {\em wins} iff $z\in A$. 
{\em Strategies} and {\em winning strategies} for one or other of the players are defined in an obvious way, and by recursive coding can be considered also elements of $\nat^{\nat}$.
$G_A$  is {\em determined} if $I$ or $I\! I$ has a winning strategy.\\

\nod Example (iii).  The {\em open game quantifier} $\game _{o}$:
$$\game _{o}uP(u) \Equi \ex x_{0}\all x_{1}\ex x_{2}\, \cdots \, \ex n P(\pa{x_{0},\ldots\, , x_{n}})$$ with dual a closed game quantifier $\game _{c}$ which we let the reader formulate. 

\section{General Infinite Game Quantifiers}

  Let  $R$ be a relation on $\nat \times \baire$. We may also use the $\game$-quantifier as embodying an operator on relations, as in the next definition.

  \begin{definition}\mb \\
$\game  \vec y R(k,\vec y)\, =_{df}
\,\,\{ k \mid \mbox{ Player $I$ has a winning strategy playing into }  \{ Y \mid R(k,Y)\}\,\}$\\
\mbox{ } \hspace{5.em}$ = \,\, \{  k \mid \exists y_{0} \all y_{1}\exists y_{2} \cdots R(k, \la y_{0},y_{i},\cdots \ra)
  \}.$
  \end{definition}

We also write this as $\game Y R(k,Y)$. For our purposes we  adapt this as follows: $\game X \Phi(k,X, \{n\mid `\psi(n,Y)\text{'}\}) $ with both $\Phi, \psi \in \call^{2}$, $Y\in \calm$, will be a new formula in $ \call^{2,\game}$
and we shall define:
\begin{center}
$\calm \models$ ``$ \game X \Phi(k,X,  `\psi(v_{0},Y)\text{'})$''$ \Longleftrightarrow  $
$ \game X \Phi(k,X, \{n\in\nat\mid (\psi(n/v_{0},Y))^{\calm}\})$.\\
\end{center}

  The latter half is to be evaluated again in $V$. There is no suggestion that $\Phi$ is absolute between the structure $\calm$ and $V$  (compare the situation with the $\mathsf{Q}^{cof}_{\omega}$ quantifier logic), nor yet that strategies (as sets of integers) are in $\calm$'s domain.  Note also that $\Phi,\psi$ are in $\call^{2}$: we are not taking $\Phi$  or $\psi$ possibly from $ \call^{2,\game}$:
there are thus (at this point) no nested $\game$ quantifiers.

{\em For $\Phi$  from a particular class $\Gamma = \Pi^{1}_{n}$ say, we are thus adding sets given by $\game \Phi$ definitions for $\Phi$ in $\Gamma$, with second order definable parameters.}

We thus define a hierarchy $\up^{\game\Gamma}_{\alpha}$ for $\alpha \in On$ using the logic $ \call^{2,\game \Gamma}$ for $\Phi\in \Gamma$:
\begin{center}
$(*)\quad \up^{\game\Gamma}_{\alpha} \models$ ``$ \game X \Phi(k,X,  `\psi(v_{0},Y)\text{'})$''$ \Longleftrightarrow  $
$ \game X \Phi(k,X, \{n\in\nat\mid (\psi(n/v_{0},Y))^{\up^{\game\Gamma}_\alpha}\})$.\\
\end{center}

By a L\"owenheim-Skolem argument, this hierarchy will close off at some countable ordinal $ \beta_{\game\Gamma}$ with resulting model
$ \up^{\game\Gamma}.$
First we consider various classes $\Gamma$.\\

{\sc Examples: (I)}  $\Gamma = \Sigma^{0}_{1}$. Then strategies for {\em open games} in  real parameters  that are $\call^{2}${\em -definable } (over $\up^{\game\Sigma^{0}_{1}}_{\alpha}$) will themselves be $ \call^{2,\game}$-definable over $\up^{\game\Sigma^{0}_{1}}_{\alpha}$ and hence will be in 
$\up^{\game\Sigma^{0}_{1}}_{\alpha+1}$. (If $I$ has a winning strategy in a  game open in the parameter $x$, then $I$ has a $HYP(x)$ winning strategy.)

 We get simply:
$$ \up^{\game\Sigma^{0}_{1}} = \up \,  \,\mbox{ and } \,\beta_{\game\Sigma^{0}_{1}}= \beta_{0}.$$
 This is because, as pointclasses,  $\game\Sigma^{0}_{1}$ $=\Pi^{1}_{1}$, and we are, in essence, defining the {\em hyperjump} (a complete $\Pi^{1}_{1}(Y)$ set of integers) over $\up_{\alpha}$ whenever $Y$ is also so definable.  $\up_{\beta_{0}}$ is closed under hyperjump, and thus for every $Y\in\power(\nat)\cap P_{\beta_{0}}$ the complete $\game\Sigma^{0}_{1}(Y)$
set is in $P_{\beta_{0}}$.  In short we are re-ordering the original Kleene hierarchy, but not adding any new sets. The same holds of the next pair of examples.\\

{\sc Examples: (II)}  $\Gamma = \Sigma^{0}_{i} \quad (i=2,3)$. 
$$ \up^{\game\Sigma^{0}_{i}} = \up \,\mbox{ and }\, \beta_{\game\Sigma^{0}_{i}}= \beta_{0} \quad\quad (i=2,3).$$
For $i=2$: it is a result of Solovay (\cf\cite{Mosch4}) that the $\game \Sigma^{0}_{2} = \Sigma^{1}_{1}$-$\mathsf{IND}$, and indeed that the least ordinal closed under such inductive definitions is less than the next $\Sigma_{2}$-admissible. Consequently Kleene's hierarchy is already closed under $\Sigma^{1}_{1}(Y)$-$\mathsf{IND}$ for any $Y\in \up$.

For $i=3$: it is a result of the author \cite{W2012} that the $\game\Sigma^{0}_{3}$ relations on integers are precisely those generalised recursive in an `eventual jump' type-2 functional $\mathsf{eJ}$, in a sense that generalises Kleene recursion in higher types. It is also shown in \cite{W2011} that the least $\Sigma_{2}$-admissible ordinal $\delta$ with $L_{\delta}$ additionally a model of $\Sigma_{2}$-Separation, is also a model of  Det({\boldmath{$\Sigma^{0}_{3} $}}). Thus $L_\delta$ is closed under 
$Y\mapsto G^{\game\Sigma^{0}_{3}}(Y)$.  And again so is $\up$.\\

In these two cases again $\up$ is the least $\beta$-model of $Z_{2}\,+$ Det({\boldmath{$\Sigma^{0}_{i} $}}) for $i=1,2,3$ so our Ramified Analytical Hierarchy has not grown by using this extra quantifier.\\

{\sc Examples: (III)}  $\Gamma = \Sigma^{0}_{i}, \,\, (i >3)$.\\  Here $ \up^{\game\Sigma^{0}_{i}} $ is still the reals of an initial segment of $L$ for some countable $ \beta_{\game\Sigma^{0}_{i}}$, but the latter ordinal is greater than $\beta_{0}$ and necessarily so by results of H. Friedman \cite{Fr70}. It is the smallest $\beta$-model of $\mathsf{Z_{2}}\, + $ Det({\boldmath{$\Sigma^{0}_{i} $}}). 

 \section{ The minimal correct model of  analysis}

Recall  (\cite{Mosch4} 6D.2)   that  the pointclass $\game \Pi^{1}_{2n-1}$ is identical to $\Sigma^{1}_{2n}$, and assuming additionally $Det(\Sigma^{1}_{2n})$, $\game \Sigma^{1}_{2n}$ is $\Pi^{1}_{2n+1}$. 
We assume from now on $PD$ or {\em Projective Determinacy} to get the right behaviour of the $\game$-quantifier on classes containing $\Pi^{1}_{1}$.

{\sc Examples: (IV)} \quad{\em The models $\up^{{ \game\Sigma^{1}_{2k}}}$ , $\up^{{ \game\Pi^{1}_{2k+1}}}$ obtained by restricting  $\game$ to be applied to formulae $\Phi \in \Sigma^{1}_{2k}$  ($\Pi^{1}_{2k+1}$ respectively) for a fixed $k$.}\\ 

Thus one obtains a model 
$\up^{{ \game\Pi^{1}_{2k+1}}}$ built \via a hierarchy of ordinal length some $\beta_{2k+1}$.  The models $\up^{{ \game\Sigma^{1}_{2k}}}$ built similarly  \via a hierarchy of ordinal length $ \beta_{2k} $ are defined analogously. 
    We shall have:
$$ 
\up^{{ \Pi^{1}_{1}}} \subset \cdots \subset \up^{{ \game\Sigma^{1}_{2k}}} \subset \up^{{ \game\Pi^{1}_{2k+1}}}\subset \cdots \subset \up^{Proj}
$$

We investigate the models $\up^{{ \game\Pi^{1}_{2k+1}}}$ a little more closely. By our assumption of $PD$ the $\Delta^{1}_{2k+2}(X)$ sets of integers from a basis for the $\Sigma^{1}_{2k+2}(X)$ relations (Moschovakis \cite{Mosch71a}). Since for any $X\in P^{{ \game\Pi^{1}_{2k+1}}}_{\alpha}$ we have any true $\Delta^{1}_{2k+2}(X)\in 
P^{{ \game\Pi^{1}_{2k+1}}}$, we may conclude that   $\up^{{ \game\Pi^{1}_{2k+1}}}$  is $\Sigma^{1}_{2k+2}$-correct. A couple of observations then follow. 
Firstly, by work of Woodin (see \cite{St95}), if
 $M^{\sharp}_{2k-1}(X)$ is the least fully iterable $X$-mouse with a measure above $2k-1$ Woodin cardinals, then it is a $\Pi^{1}_{2k+1}$-singleton set, and thus has a code as a $\Delta^{1}_{2k+2}$ set of integers (such mice exist thanks to $PD$). Thus $\up^{{ \game\Pi^{1}_{2k+1}}}$ is closed also under $X\imp M^{\sharp}_{2k-1}(X)$. 
 Secondly, recall that:
 $$\delta^{1}_{n}= \sup \{\tmop{rk}( R ) \mid R \mbox{ a } \Delta^{1}_{n} \text{-prewellording of } \nat \}.$$ 
 (Recall also that $\delta^{1}_{1} =\omega_{1}^{ck}$ and $\delta^{1}_{2} =\sigma_{1}$ the least $\Sigma_{1}$ stable ordinal.)
Thus $\up^{{ \game\Pi^{1}_{2k+1}}}$ is also closed under $X\imp \delta^{1}_{2n+2}(X)$. We should probably point out that in no real sense is the hierarchy  $\up^{{ \game\Pi^{1}_{2k+3}}}$ an end-extension of $\up^{{ \game\Pi^{1}_{2k+1}}}$: 
the sets in the latter appear all at the first or second stage of the former.
\\

The minimal $\Sigma^{1}_{2n}$-correct models were first identified by Enderton and H. Friedman \cite{EndFr71}. They built their models, for a given $n$, and obtained their correctness, by {\em assuming} that $\Sigma^{1}_{{2n}}$ relations had a {\em basis} in the $\Delta^{1}_{2n}$ definable reals (that is every $\Sigma^{1}_{{2n}}$ relation contains a $\Delta^{1}_{2n}$ definable point). In 1971, as they noted in their conclusion, it was still a conjecture that $PD$ implied this latter basis result. They also noted a conjecture of Martin and Solovay, which also turned out to be true under $PD$: that  $\Sigma^{1}_{{2n+1}}$ relations did not have a {basis} in the $\Delta^{1}_{2n+1}$ definable reals. (The correct statement, under $PD$, is that the set of reals recursive in the real of some (equivalently, of any)  $\Pi^{1}_{2n+1}$, but not $\Delta^{1}_{2n+1}$, definable singleton set, form a basis for $\Sigma^{1}_{{2n+1}}$, \cf \cite{KeMaSo88} or \cite{Mosch4} 6C.10.) Thus under PD $\up^{{ \game\Pi^{1}_{2k+1}}}$ is $\Sigma^{1}_{2k+2}$-correct. But the failure of the unamended basis theorem implies that $\up^{{ \game\Sigma^{1}_{2k+2}}}$ is not $\Sigma^{1}_{2k+3}$-correct.

(They also performed their construction of $\Sigma^{1}_{2n+3}$-correct minimal models, whilst hypothesizing the (false under $PD$) basis assumption that 
$\Sigma^{1}_{{2n+1}}$ relations did have a {basis} in the $\Delta^{1}_{2n+1}$ definable reals. They remarked that the hypothesis is after all {\em consistent} since it holds in $L[\mu]$  - but seemed not to notice that in fact it holds in any case in $L$. As a final historical remark Shilleto \cite{Sh72} constructed in a slightly complicated fashion minimal  $\Sigma^{1}_{n}$-correct models but assumed $V=L$.  For $n=2$ the Enderton-Friedman construction is simpler.)

\begin{theorem}
$\up^{{ \game\Pi^{1}_{2k+1}}}$ is the minimal $\Sigma^{1}_{2k+2}$-correct model of analysis.
\end{theorem}
Proof: We have seen correctness above. The issue is only minimality. We only outline the steps.
 Let $\calm$ be any other such $\Sigma^{1}_{2k+2}$-correct model of analysis. Let $\beta_{2k+1}$ be the closure ordinal of the hierarchy of the model $ \up:=\up^{{ \game\Pi^{1}_{2k+1}}}$. (We abbreviate for this proof. $ \up^{{ \game\Pi^{1}_{2k+1}}}_{\alpha}$ as $\up_{\alpha}$.)  The idea is just that of defining $\la \up^{M}_{\alpha} \mid \alpha \leq W_{\calm} \ra$ where $W_{\calm}$ is the supremum of ordinals representable by reals of $\calm$. One thus should show that we can define within second order arithmetic the $P_{\alpha}$ hierarchy within \calm  for any $\alpha$ representable in \calm and in an absolute fashion which ensures $\up_{\alpha}=\up_{\alpha}^{\calm}$. The papers of \cite{BHP66} and \cite{EndFr71} give in great detail how this may be done in the simpler second order number-theoretic sense (in the first paper), and using additional assumptions of basis theorems (in the second paper). We shall assume that readers will believe that such formalisations are possible without wearing them out with the details here. Perhaps there are two points to be emphasised here. The first is in the transition from 
$P_{\alpha}$ to $P_{\alpha+1}$ that by the $\Sigma^{1}_{2k+2}$-correctness of $\calm$ (and Moschovakis' Third Periodicity Theorem) all the necessary strategies needed to give the correct evaluation of a formula are available in the model \calm . 
(In more detail: Third Periodicity says that any $\Pi^{1}_{2k+1}(X)$ game that is won by Player $I$ has a winning strategy that is $\Sigma^{1}_{2k+2}(X)$. But as a consequence that winning strategy will be in $\up^{{ \game\Pi^{1}_{2k+1}}} $ if $X$ is.)
Thus we shall have (for $\Phi\in \Gamma =\Pi^{1}_{2k+1}$ and assuming inductively that ${\up_\alpha} = ({\up_\alpha})^{\calm}$):\\

$ \game X \Phi(k,X, \{n\in\nat\mid (\psi(n/v_{0},Y))^{\up^{\game\Gamma}_\alpha}\})  \Longleftrightarrow  $\\
\mb \hfill $ \Longleftrightarrow ( \game X \Phi(k,X, \{n\in\nat\mid (\psi(n/v_{0},Y))^{\up^{\game\Gamma}_\alpha}\}))^{\calm}$\\

and hence $P_{\alpha+1} = P_{\alpha+1}^{\calm}$. We shall not say any more on this point.

We may set $\bar \beta \dfs W_{\calm}$. Clearly $Lim(\bar\beta)$. Then $\beta_{2k+1}\leq\bar\beta$ must hold thus establishing the required minimality as then $\up$ is an initial segment of $ \up^{\calm}$. For if $\beta_{2k+1}>\bar\beta$, we should have a failure of comprehension (in our expanded sense, meaning closure of definability in the extended logic) over $\up_{\bar\beta}$ and thence over $ \up_{\bar\beta}^{\calm} $. The latter is a definable {\em `class'} of $\calm$,  not being coded by any set of $\calm$.
Indeed that failure of comprehension can be strengthened to show that we actually have a code for a wellordering $u$ of type $\bar \beta$ definable over $\up_{\bar\beta}$. However then we have a wellorder $u$ of order type $\bar \beta$ which is definable over $\up_{\bar\beta}=\up^{\calm}$. But the latter is a definable class of $\calm$; so $u$ is definable over $\calm$, and so must be in \calm as the latter is a model of full $\Pi^{1}_{\omega}$-comprehension, thus leading to an obvious contradiction.
\hfill \qed \\

The last argument shows immediately:
\begin{corollary}\label{4.2} $\beta_{2k+1}= sup \{ rk(Y)\mid Y\in WO\cap\up^{{ \game\Pi^{1}_{2k+1}}}\}$.
\end{corollary} 
{\sc Examples: (V)} {\em Let } $\Gamma = \Pi^{1}_{\omega} = Proj.$\\
So now $\Df_{\call^{2,\game}}$ yields a model $$\up^{Proj} = \up^{\game\Pi^{1}_{\omega}}. $$
Let  $\beta_{Proj} $ be the closure ordinal of $\up^{Proj}$.

\begin{theorem}[PD]
$ \up :=\up^{Proj}$ is the minimal projectively correct model of Analysis (and so also the minimal projectively correct model of Analysis + {{$PD$}}). Moreover 
$$\up=\bigcup_{k} \up^{\game\Pi^{1}_{2k+1}}
 = \bigcup_{k} \up^{\game\Sigma^{1}_{2k}}$$
\end{theorem}

\rem By ``$PD$'' we mean the scheme that contains for every $n\in\nat$ the statement that  {\em ``For every $\Pi^{1}_{n}$ set $A\sset \baire$, and tree $T\sset \omega^{<\omega}$,  $G(A,T)$ is determined.''}

\nod{\sc Proof:} 
We have for each $k$: 
$ \up^{\game\Pi^{1}_{2k-1}}$  is the minimal $\Sigma^{1}_{2k}$-correct model. We may naturally write:
$ \up^{\game\Pi^{1}_{2k-1}}\prec _{\Sigma^{1}_{2k}}  \cal{Z}$. From this it follows easily that 
$$\bigcup_{k} \up^{\game\Pi^{1}_{2k-1}}\prec _{\Sigma^{1}_{\omega}}\cal{Z}$$ is the minimal fully correct model of analysis. The same follows for $\bigcup_{k} \up^{\game\Sigma^{1}_{2k}}$. By definition $\up^{Proj} \supset \up^{\game\Sigma^{1}_{2k}},  \up^{\game\Pi^{1}_{2k+1}}$. An induction on $\alpha <\beta_{Proj}$ shows that $\up^{Proj}_{\alpha}\sset \bigcup_{k} \up^{\game\Pi^{1}_{2k+1}}$ also.  The statement ``{\em Determinacy}$(${\boldmath$\Pi$}$^{1}_{n})$'' is expressible by a projective formula, and as true in $\cal{Z}$ it will be true in $\up^{Proj}$.
\hfill\qed\\

To give some further idea of what these models contain we use further descriptive set theoretical ideas. 
\begin{theorem} [Woodin]{\em(\cite{St95})} ($\all n \,M_n^\sharp$ {\em exists}) \quad $\re\cap M_{2n} = \calc_{2n}$; \quad $\re \cap M_{2n-1} = Q_{2n+1}$.
\end{theorem}
Here $C_{2n}$ is the largest countable $\Sigma^{1}_{2n}$ set of reals;  $Q_{2n+1}$ is the set of reals each of which is $\Delta^{1}_{2n+1}$ definable in (a code for) a countable ordinal. Because $\up^{Proj}$ contains for every $X$ (a code for) $M_n^{\sharp}(X) $, it will in particular contain all the reals of $M_{n}(X)$ for each $n$ (as such reals are all recursive in $M_n^{\sharp}(X) $); relativising the last result we shall have:

\begin{corollary}
For every $n$, 
 for all $X \in P^{Proj}$: $$\calc_{2n}(X), Q_{2n+1}(X)\sset P^{Proj}.$$
\end{corollary}

\begin{corollary} $\up:=\up^{\game\Pi^{1}_{2k+1}}$ is closed under
$$ \all X \in P, \,\calc_{2k-2}(X), Q_{2k+1}(X)\sset P.$$ 
\end{corollary}

We now try to identify the reals of $\up^{Proj}$ in terms of a level of an inner model (just as for Gandy-Putnam, the original ramified hierarchy continued for $\beta_{0}$ steps, and whose reals were precisely those of $L_{\beta_0}$). Let $\beta_{{Proj}}$ be the closure ordinal of the $P^{Proj}_{\alpha}$ hierarchy. We use another result of Woodin:

\begin{theorem}[Woodin]{\em (\cite{St95}, 4.7)}\label{4.7} Assume $M_{2n}^{\sharp}$ exists. Then $M_{2n}$ is $\Sigma^{1}_{2n+2}$-correct.
\end{theorem}

Hence  the  $\up^{\game\Pi^{1}_{2k+1}}_{\alpha}$ hierarchy will be absolute between $V$ and $M_{2k}$. Hence:
$$\mbox{(a)  } \quad\up^{\game\Pi^{1}_{2k+1}}_{\alpha}\sset  M_{2k}^{\sharp}.$$
Note also that 
$$\mbox{(b) }\quad M_{2n}^{\sharp}\notin \up^{\game\Pi^{1}_{2k+1}} \mbox{ but it is in } \up^{\game\Pi^{1}_{2k+3}}.$$

Let $M_{\omega}= L[E_{\omega}]$ be the minimal iterable model of $\omega$ Woodin cardinals. (This is somewhat overkill: we need simply a hierarchy containing all the $M_{n}^{\sharp}$'s.) Let $\tau_{n}$ be the index in $L[E_{\omega}]$ which attaches the topmost sharp filter for the model $M_{2n}$. In other words, so that $M_{2n}^{\sharp}= \la J_{\tau_{n}}^{E_{\omega}},\in, E^{\omega}\!\restriction \!\tau_{n}, F_{\tau_{n}}\ra$ (where we have followed the usage of the Jensen $J$-hierarchy when defining such models). Let $\tau= \sup_{n}\tau_{n}$. We let $Q\dfs \la J_{\tau}^{E_{\omega}},\in\ra$ be union of these levels. Then for any $k$, $Q$ is closed under $X\imp M_{k}^{\sharp}(X)$.

\begin{theorem} 
  $P^{Proj}=\re\cap Q$.
\end{theorem}   

\nod Proof:
Note that $Q\models$ ``$V=HC$'' (as there is definably over each structure $\la J_{\tau_{n}}^{E_{\omega}},\in, E_{\omega}\!\restriction\! \tau_{n}, F_{\tau_{n}}\ra$ an
onto map from $\omega$  onto its domain).
As $\up^{\game\Pi^{1}_{2k+1}}_{\alpha}\sset \re \cap  J_{\tau_{n}}^{E_{\omega}} $, (by (a) following on Theorem \ref{4.7}) we have that $\up^{Proj}\sset \re \cap Q$.  Conversely any real $Y\in  \re \cap Q$ is in some 
$M_{2n}^{\sharp}$ and the latter is in $\up^{\game\Pi^{1}_{2k+3}}$ by (b). 
Hence: $\up^{Proj}=\re \cap Q$. \mb \hfill \qed

 \section{ Some intermediate models}

To discuss further models we introduce Moschovakis's notion of a {\em Spector class} of pointsets in $\power(\nat)$. Broadly speaking this
is a notion of a class of sets of integers arising from some general abstract notion of definability. Such a family must exhibit a number of properties: (i) Some elementary closure; (ii) be $\omega$-parametrized; (iii) have the all-important  {\em Scale Properties}.  The reader is referred to \cite{Mosch4} for a full definition and discussion. The following are all examples of this notion, starting with the least, and canonical, one:
$$\Pi^{1}_{1}= \game\Sigma^{0}_{1}\,;\, \game\Sigma^{0}_{n}\,; \,\Sigma^{1}_{2}\,;  \,\mbox{(and under $PD$) }\, \Pi^{1}_{2n+1}, \,\Sigma^{1}_{2n+2} \ldots$$
Initially we disbarred nested applications of the game quantifier to formulae. We can lift that restriction to obtain another sub-hierarchy of models.

For $\Gamma'$ a Spector class we set $\Delta' = \Gamma' \cap \widecheck{\Gamma'}$ to be  the {\em self-dual} part of $\Gamma'$. For Spector classes $\Gamma, \Gamma'$ we set $\Gamma \prec \Gamma'$ iff 
$ \Gamma\sset{\Delta'}.$

\begin{definition}(Spector Ordinal)
Let $\Gamma$ be a Spector Class, $\Delta = \Gamma \cap \widecheck{\Gamma}$ its self-dual part. We set
$$
\kappa^{\Gamma}\dfs \sup \{ \tmop{rk}(P)\mid P\in \Delta, P \mbox{ a prewellordering of } \omega \}.
$$

\end{definition}

\begin{lemma}[Moschovakis](Spector Criterion)
Let $\Gamma, \Gamma'$ be two Spector classes on $\nat$.
$$
\Gamma\sset \Gamma' \imp (\Gamma \prec \Gamma' \,\equi \,\kappa^{\Gamma} <\kappa^{\Gamma'})
.$$
\end{lemma}

\nod {\sc Examples: (VI)} \quad  {\em Allow formulae with nested $\game$ quantifiers.}\\ Let $\game^{n}\Gamma$ be the pointclass of sets defined by formulae of the form $\game\cdots\game\Phi$ for a $\Phi \in \Gamma$. These are also Spector pointclasses. Set $\Gamma_{2n+1}=\Pi^{1}_{2n+1}$ and $\Gamma_{2n}= \Sigma^{1}_{2n}$. We adopt the abbreviation:  $\Gamma _{k,n}=
\game^{n}\Gamma_{k}$. Then each  $\Gamma _{k,n}$ is a Spector pointclass and
$$
 \Pi^{1}_{2k+1}= \Gamma_{2k,1} \prec  \Gamma_{2k,2}\prec   \cdots
\prec \Gamma_{2k,n} \prec \cdots  \Sigma^{1}_{2k+2}
 $$
 \nod with corresponding models:
 $$
  \up^{\Gamma_{2k,1}} \subseteq \up^{\Gamma_{2k,2}}\subseteq \cdots 
\subseteq   \up^{\Gamma_{2k,n}}\subseteq \cdots
 \subseteq \up^{\Pi^{1}_{2k+3}}
$$
and similarly for $\Gamma_{2k+1}$ {\em m.m.} (see \cite{Ke78} 2.5.2).

 \section{ Generalised logics from Spector Classes: some conclusions}

Seeing that game quantifiers give rise to Spector classes when applied to Spector classes, one may start to query whether a given notion of definability as encapsulated in an abstract Spector class  $\Gamma$ can be in turn  used to create a model of analysis, $\up^{\Gamma}$, in this way. The thought is that the notion of generalised definability can be applied in a hierarchical fashion, level by level, to bring more sets into the hierarchy. But how is this to be done, or rather, can it be done in a manner that fits this investigation of extended logics coupled with a ramified approach? The answer turns out that this can indeed be done and relatively easily, thanks to the following result.

\begin{theorem}[Harrington (\cite{Ke78} 3.2)]
 Let ${\Gamma}\sset \power(\nat)$ be a Spector Class. Then there is a generalized quantifier  $\mathsf{Q}$ so that ${\Gamma} = IND(\mathsf{Q})$.
\end{theorem}

By  $IND(\mathsf{Q})$ we mean the class of sets of integers inductive using now formulae in $\call^{1,\mathsf{Q}}$ over $\la  \nat, +,\times, \ldots \ra$. More specifically we adjoin to the first order language the quantifier $\mathsf{Q}$, call this $\call^{1,\mathsf{Q}}$, and consider an $\call^{1,\mathsf{Q}}$ formula $\varphi(v_{0},S)$ where $S$ is a second order variable which only appears positively in $\varphi$ - that is within an even number of negations.  (See Moschovakis \cite{M} Ch.9.) Then one may build up successive extensions of $S_{0}=\emp$, $S_{\alpha+1}= \{ n\mid  \la  \nat, +,\times, \ldots , S_{\alpha}\ra \models \varphi [n,S_{\alpha}]\}$ in the familiar fashion.  By the positivity requirement on $S$ this is a monotone increasing hierarchy, which by taking unions at limits, $S_{\lambda}$, reaches a fixed point $S_{\infty}$. The theorem above then says that for any Spector Class $\Gamma$ there is a corresponding $\mathsf{Q}_
{\Gamma}$ which will inductively define in $\call^{1,\mathsf{Q}_{\Gamma}}$ precisely all and only the members of $\Gamma$.

\begin{theorem}\label{6.2} Let $\Gamma\sset \power(\nat)$ be a Spector class, with corresponding quantifier $\mathsf{Q}= \mathsf{Q}_{\Gamma}$ from the last theorem. Then there is a minimum model of analysis $\up^{\Gamma}$ which is closed under positive inductions in $\call^{\mathsf{Q}}$, and so that for any $X\in P^{\Gamma}$ we have $\Gamma(X)$ (the Spector class relativised to $X$) is contained in $\up^{\Gamma}$.
\end{theorem}
Proof: It should be clear that for any parameter $X\in\power(\nat)\cap P^{\Gamma}_{\alpha}$ that if $S_{\beta}$ is a stage in some inductive definition using some $\call^{1,\mathsf{Q}}$ formula $\varphi(v_{0},X,S)$ with $S_{\beta}\in \up_{\alpha}$, then $S_{\beta+1}$ will be placed in $ \up^{\Gamma}_{\alpha+1}$. Thus $P^{\Gamma}_\alpha$ continues to grow until we  reach a closure point that contains all fixed points for all such inductions using all possible parameters. \hfill \qed\\

Moral: we can add an abstract quantifier  $\mathsf{Q_{\Gamma}}$ to obtain a language $\call^{2,\mathsf{Q_{\Gamma}}}$ to close up under inductions in $\mathsf{Q_{\Gamma}}$ and so define a ramified hierarchy  using the kind of definability given by $\Gamma$. In other words, Spector classes $\Gamma$ give rise to models $\up^{\Gamma}$ which are minimum models of analysis closed under $X\imp\,\Gamma(X)$. It is in this sense that we claim, as in our introduction, that if we identify plausible notions of definability with Spector classes, then we can find a quantifier, and hence an extended logic, to build a ramified hierarchy to exemplify it. \\

We conclude with some loose ends in the guise of open questions:\\

{\em Question: Characterise for which Spector classes $\Gamma,  \Gamma'$ we have that $$\Gamma\prec  \Gamma' \,\Imp\, \up^{\Gamma}\subset \up^{\Gamma'}.$$} 
\nod Note that this is a non-trivial question: for $\Gamma_{i} = \game\Sigma^{0}_{i}$ $(0<i<4)$, we have that $\Gamma_{i}\prec \Gamma_{i+1}$ but all three classes $\Gamma_{i}$ have the same $\up^{\Gamma_{i}}$ namely Kleene's original $\up$.
Here is another specific case where we do not know the answer:\\

{\em Question: For a fixed $k$, are the inclusions between the models $\up^{\Gamma_{2k,n}}$ in Examples (VII) strict?} \\

We have left open the question of models with an odd levels of correctness in the projective hierarchy:\\

{\em Question: Assume $PD$. Identify a logic which builds  the  minimal $\Sigma^{1}_{2k+1}$-correct model of analysis.}\\

{\em Question: Assume $PD$. Is there some further characterization of the length of the ordinals $\beta_{2k},\beta_{{2k+1}}$?}
\small
\bibliographystyle{plain}
\bibliography{settheory10s}

\ed